\documentclass[a4paper,11pt]{article}
\usepackage{amssymb}
\usepackage{amsmath,amsfonts,amsthm,amssymb}
\usepackage[dvips]{graphics}
\usepackage{epsfig}
\usepackage{indentfirst}
\usepackage{color}
\allowdisplaybreaks

\newtheoremstyle{theorem}
  {10pt}          
  {10pt}  
  {\sl}  
 {}
  {\bf}  
  {. }    
  { }    
  {}     
\theoremstyle{theorem}

\newtheorem{theorem}{Theorem}[section]

\newtheorem{definition}{Definition}[section]
 \newtheorem{lemma}{Lemma}[section]
 
 \newtheorem{remark}{Remark}[section]
  
\numberwithin{equation}{section}

\newtheoremstyle{defi}
  {10pt}          
  {10pt}  
  {\rm}  
  {}  
  {\bf}  
  {. }    
  { }    
  {}     
\theoremstyle{defi}

\begin{document}
\baselineskip = 16pt

\title{\bf On the hydrostatic approximation of compressible anisotropic Navier-Stokes equations - rigorous justification }

\author{ Hongjun Gao$^{1}$  \footnote{Email:gaohj@hotmail.com}\ \ \
\v{S}\'{a}rka Ne\v{c}asov\'{a}$^2$
\footnote{Email: matus@math.cas.cz} \ \ \  Tong Tang$^{3}$ \footnote{Email: tt0507010156@126.com}\\
{\small 1.School of Mathematical Sciences,}\\
{\small Southeast University, Nanjing 210096, P.R. China}\\
{\small  2. Institute of Mathematics of the Academy of Sciences of the Czech Republic,} \\
{\small \v Zitn\' a 25, 11567, Praha 1, Czech Republic}\\
{\small 3. Department of Mathematics, College of Sciences,}\\
{\small Hohai University, Nanjing 210098, P.R. China}\\
\date{}}

\maketitle
\begin{abstract}
In this work, we derive the hydrostatic approximation by taking the small aspect ratio limit to the Navier-Stokes equations. The aspect ratio (the ratio of the depth to horizontal width) is a geometrical constraint in general large scale geophysical motions that the vertical scale is significantly smaller than horizontal. We derive the versatile relative entropy inequality. Applying the versatile relative entropy inequality we gave the  rigorous derivation of  the limit from the compressible Navier-Stokes equations to the compressible Primitive Equations. This is the first work where the relative entropy inequality was used for proving hydrostatic approximation - the compressible Primitive Equations.
\vspace{0.5cm}


\medskip

\end{abstract}

\maketitle
\section{Introduction}\setcounter{equation}{0}
{The atmosphere and ocean have attracted considerable attention in the scientific research community, especially for the geophysics, as it has so many fluid dynamic properties and mysterious phenomena. One of the most interesting and physically important features of large-scale meteorology and oceanography is  that vertical dimension of the domain is much smaller than the horizontal dimension of domain.} Therefore, many scientists have suggested the viscosity coefficients must be anisotropic, such as \cite{bry,pe,was}. The anisotropic Navier-Stokes equations are widely used in geophysical fluid dynamics. In this paper, we consider the following compressible anisotropic Navier-Stokes equations
\begin{eqnarray}
\left\{
\begin{array}{llll}  \partial_{t}\rho+\text{div}(\rho \mathbf{u})=0,\\
(\rho\mathbf u)_t+\text{div}(\rho\mathbf u\otimes\mathbf u)+\nabla p(\rho)=\mu_x\Delta_x\mathbf u+\mu_z\partial_{zz}\mathbf u,
\end{array}\right.
\end{eqnarray}
in the thin domain $(0,T)\times\Omega_\epsilon$. Here $\Omega_\epsilon=\{(x,z)|x\in\mathbb{T}^2,-\epsilon<z<\epsilon\}$, $x$ denotes the horizontal direction and $z$ denotes the vertical direction, while, $\mu_x$ and $\mu_z$ are given constant horizontal viscous coefficient and vertical viscous coefficient. The velocity $\mathbf{u}=(\mathbf v,w)$, where $\mathbf v(t,x,z)\in\mathbb{R}^2$ and $w(t,x,z)\in\mathbb{R}$ represent the horizonal velocity and vertical velocity respectively. Through out this paper, we use $\rm div\mathbf u=div_{x}\mathbf v+\partial_zw$ and $\nabla=(\nabla_x,\partial_z)$ to denote the three-dimensional spatial divergence and gradient respectively, and $\Delta_x$ stands for horizontal Laplacian. As atmosphere and ocean are the thin layers, where the fluid layer depth is small compared to radius of sphere, Pedlosky \cite{pe} pointed out that "the pressure difference between any two points on the same vertical line depends only on the weight of the fluid between these points...". Here we neglect the gravity and suppose the pressure $p(\rho)$ satisfies the barotropic pressure law where the pressure and the density are related by the formula: $p(\rho)=\rho^\gamma\hspace{5pt}(\gamma>1)$. {\bf Therefore we assume the density $\rho$ is independent of $z$} that is $\rho=\rho(t,x)$. This plausible assumption agrees well with experiment and is frequently taken as a hypothesis in geophysical fluid dynamics.

Similar to the assumptions by \cite{a,li}, we suppose $\mu_x=1$ and $\mu_z=\epsilon^2$. As stressed by Az\'{e}rad and Guill\'{e}n \cite{a}, it is necessary to consider the above anisotropic viscosities scaling, which is fundamental for the derivation of Primitive Equations (PE).{\footnote{For completeness we will write the system (PE):
\begin{eqnarray}
\left\{
\begin{array}{l}
\partial_t \mathbf{v}+\textrm{div}_x(\mathbf{v}\otimes\mathbf{v})+\partial_z (\mathbf vw)+\nabla_x p=\Delta_x\mathbf v+\partial_{zz}\mathbf v,\\
\partial_zp=0\\
\text{div}_x \mathbf{v}+\partial_z w=0,
\end{array}\right.
\end{eqnarray}}}

 Under this assumption, the system is rewritten as the following
\begin{eqnarray}
\left\{
\begin{array}{llll}  \partial_{t}\rho+\text{div}_x(\rho \mathbf{v})+\partial_z(\rho w)=0, \\
\rho\partial_t\mathbf{v}+\rho(\mathbf u\cdot\nabla)\mathbf v-\Delta_x\mathbf v-\epsilon^2\partial_{zz}\mathbf v+\nabla_x p(\rho)=0,\\
\rho\partial_tw+\rho\mathbf u\cdot\nabla w-\Delta_x\mathbf w-\epsilon^2\partial_{zz}w+\partial_zp(\rho)=0.
\end{array}\right.\label{a}
\end{eqnarray}

Inspired by \cite{li}, we introduce the following new unknowns,
\begin{eqnarray*}
\mathbf u_\epsilon=(\mathbf v_\epsilon,w_\epsilon),\hspace{2pt}\mathbf v_\epsilon(x,z,t)=\mathbf v(x,\epsilon z,t),\hspace{2pt}
w_\epsilon=\frac{1}{\epsilon}w(x,\epsilon z,t),\hspace{2pt}\rho_\epsilon=\rho(x,t),
\end{eqnarray*}
for any $(x,z)\in\Omega:=\mathbb{T}^2\times(-1,1)$. Then the system (\ref{a}) becomes the following compressible scaled Navier-Stokes equations (CNS):
\begin{eqnarray}
\left\{
\begin{array}{llll}  \partial_{t}\rho_\epsilon+\text{div}_x(\rho_\epsilon\mathbf{v}_\epsilon)+\partial_z(\rho_\epsilon w_\epsilon)=0, \\
\rho_\epsilon\partial_t\mathbf{v}_\epsilon+\rho_\epsilon(\mathbf u_\epsilon\cdot\nabla)\mathbf v_\epsilon-\Delta_x\mathbf v_\epsilon-\partial_{zz}\mathbf v_\epsilon+\nabla_x p(\rho_\epsilon)=0,\\
\epsilon^2(\rho_\epsilon\partial_tw_\epsilon+\rho_\epsilon\mathbf u_\epsilon\cdot\nabla w_\epsilon-\Delta_x\mathbf w_\epsilon-\partial_{zz}w_\epsilon)+\partial_zp(\rho_\epsilon)=0.
\end{array}\right.\label{aa}
\end{eqnarray}
We supplement the CNS with the following boundary and initial conditions:
\begin{eqnarray}
&\rho_\epsilon,\hspace{2pt}\mathbf u_\epsilon\hspace{2pt}\text{are periodic in}\hspace{2pt} x, y, z,\nonumber\\
&(\rho_\epsilon,\mathbf u_\epsilon)|_{t=0}=(\rho_0,\mathbf u_0).\label{pa}
\end{eqnarray}
The goal of this work is to investigate the limit process $\epsilon\rightarrow0$ in the system of (\ref{aa}) converge in a certain sense to the following compressible Primitive Equations(CPE):
\begin{eqnarray}
\left\{
\begin{array}{llll}  \partial_{t}\rho+\text{div}_x(\rho \mathbf{v})+\partial_z(\rho w)=0, \\
\partial_t(\rho \mathbf{v})+\textrm{div}_x(\rho\mathbf{v}\otimes\mathbf{v})+\partial_z(\rho\mathbf vw)+\nabla_x p(\rho)=\Delta_x\mathbf v+\partial_{zz}\mathbf v,\\
\partial_zp(\rho)=0.
\end{array}\right.\label{aaa}
\end{eqnarray}

The geophysical fluid dynamics is a fundamental subject to understand the atmosphere and ocean. Whereas, from the mathematical point of view and numerical perspective, it is very complicated to use the full hydrodynamical and thermodynamical equations to analyze and simulate  { atmospheric} flows and oceanic flows. Therefore, scientists have  introduced the Primitive Equation (PE) model in the geophysical fluid dynamics.{\footnote{By (PE) we mean the {\it incompressible} Primitive Equation.}} It was Richardson that derived originally PE model in 1920's for weather prediction. But lacking stability of calculations, this model was not so successful. Then, Bryan \cite{bry} improved PE model by applying the hydrostatic approximation in 1969. Compared with abundant successful results in simulation and application for PE at early stage, the mathematical research of PE was started very late. It was until 1990s that Lions, Teman and Wang \cite{l1,l2} were first to study the PE and received fundamental results in this field. Then PE has historically progressed by concentrated the mathematical arguments developed by the precise analysis of simpler models. There is a lot of works  dedicated to PE model see \cite{bg,b1,b2,c2,c4,c5,g2,g3,ju,l3,l5,t,ws} and references therein. Let us give a short retrospect and comment for some results. Guill\'{e}n-Gonz\'{a}lez, Masmoudi and Rodr\'{\i}guez-Bellido \cite{gu} proved the local existence of strong solutions in the three dimension case. The celebrated breakthrough result was made by Cao and Titi \cite{c1}. They were the first who proved the global well-posedness of PE in the three dimensional case. Then, by virtue of semigroup method, Hieber and Kashiwabara \cite{h} extended this result relaxing the smoothness on the initial data. On the other hand, regarding to inviscid PE (hydrostatic incompressible Euler equations), the existence and uniqueness is an outstanding open problem. Brenier \cite{b} proved the existence of smooth solutions in two-dimensions under the convex horizontal velocity assumptions. And he in \cite{by} suggested that the existence problem may be ill-posed in Sobolev spaces. Later, Masmoudi and Wong \cite{m} extended Brenier's result, removing the convex horizontal velocity assumptions. Partly for historical reasons, the research of geophysical fluid concerns on PE model at incompressible case. However, it is well known that atmosphere and ocean have compressible property. Therefore, it is interesting and natural to consider the PE model at compressible case, that is CPE. With the constant viscosity coefficients, Gatapov and Kazhikhov \cite{g}, Ersoy and Ngom \cite{er2} proved the global existence of weak solutions in 2D case. Recently, Liu and Titi \cite{liu1,liu3} proved the local existence of strong solutions in 3D case and consider the zero Mach number limit of CPE. On the other hand, Ersoy et al. \cite{er1} used the dimensionless number and asymptotic analysis, obtaining the CPE in the case where the viscosity coefficients are depending on the density. Ersoy et al. \cite{er1}, Tang and Gao \cite{tang} showed the stability of weak solutions. The stability means that a subsequence of weak solutions will converge to another weak solutions if it satisfies some uniform bounds. Recently, Liu and Titi \cite{liu2} and independently  Wang et al. \cite{w} used the B-D entropy to prove the global existence of weak solutions.

As stressed by \cite{a,li}, the hydrostatic approximation is one of the important feature of PE model. A rigorous justification of the limit passage from anisotropic Navier-Stokes equations to its hydrostatic approximation via the small aspect limit seems to be of obvious practical importance. There are numerous studies of the incompressible convergence. For example, Az\'{e}ard and Guill\'{e}n \cite{a} proved the weak solutions of anisotropic Navier-Stokes converges to weak solutions of PE. Li and Titi \cite{li} used the method of weak-strong uniqueness to prove the aspect ratio limit of incompressible anisotropic Navier-Stokes equations, that is from weak solutions of anisotropic Navier-Stokes equations to strong solutions of incompressible PE model. Then Giga, Hieber and Kashiwabara et al. \cite{f1,f2} extended the results into maximal regularity spaces. Recently, Donatelli and Juhasz \cite{do} proved the convergence in downwind-matching coordinates. For the stationary case, readers can refer to \cite{be,b2}. On the other hand, based on a revised global Cauchy-Kowalewski theorem, Paicu, Zhang and Zhang \cite{p} proved the incompressible anisotropic Navier-Stokes equations converge to the Prandtl equation in Besov spaces for 2D case. However, for the compressible fluids flows, {\bf to the best of authors' knowledge, there are no results concerning the convergence from compressible Navier-Stokes system (CNS) to compressible Primitive Equations (CPE)}.

Our goal is to rigorously justify the limit in the framework of weak solutions of CNS. Recently, Bella, Feireisl and Novotn\'{y} \cite{bel}, Maltese and Novotn\'{y} \cite{ma} proved the limit passage from 3D compressible Navier-Stokes equations to 1D and 2D compressible Navier-Stokes equations in thin domain. See also result by Ducomet et al. \cite{du}. Heuristically, inspired by their works, we develop and adapt the corresponding idea of relative entropy inequality for compressible Navier-Stokes equations. {\bf There are significant  differences at the mathematical structure between Navier-Stokes equations and CPE model.} Due to the hydrostatic approximation, {\bf there is no information for the vertical velocity in the momentum equation of CPE model, and the vertical velocity is determined by the horizontal velocity via the continuity equation}, so it is very difficult to analyze the CPE model. Therefore, the classical method used in Navier-Stokes system can not be applied straightforwardly to CPE. Luckily, based on our previous work \cite{gao} of weak-strong uniqueness to CPE, we prove the aspect ratio limit of compressible anisotropic Navier-Stokes equations. This is the first work to use the relative entropy inequality for proving the hydrostatic approximation at the compressible case. For the introduction of  the versatile relative entropy inequality, see \cite{gao1}. And we will introduce the differences between our present work with \cite{gao} in Section 3 Remark 3.3. Last but not least, let us mention that the corner-stone analysis of our results is based on the relative energy inequality which was invented by Dafermos, see \cite{D}. Then Germain \cite{ge} introduced it into compressible Navier-Stokes equations. It is Feireisl and his co-authors \cite{e2,e3,e4} that generalized the relative energy inequality for solving various compressible fluid model problems.

The paper is organized as follows. In Section 2, we recall some useful inequalities. We introduce the definition of weak solutions, strong solution, relative energy and state  { the main} theorem in Section 3. Section 4 is devoted to proof of the convergence.
\vskip 0.5cm
\section{Preliminaries}

In this section, we first introduce some basic inequalities needed in the later proof. The first inequality is the so called the generalized Poincar\' e inequality.
\begin{lemma}
Let $2\leq p\leq6$, and $\rho\geq0$ such that $0<M\leq\int_\Omega\rho dx,\int_\Omega\rho^\gamma dx\leq E_0$ for some $(\gamma>1)$ then
\begin{align*}
\|f\|_{L^p(\Omega)}\leq C\|\nabla f\|_{L^2(\Omega)}+\|\rho^{\frac{1}{2}}f\|_{L^2(\Omega)},
\end{align*}
where $C$ depends on $M$ and $E_0$.
\end{lemma}
The details of proof can be seen at Feireisl's monograph \cite{e1}.  The following is the famous Gagliardo-Nirenberg inequality.
\begin{lemma}
For a function $u:\Omega\rightarrow\mathbb{R}$ defined on a bounded Lipschitz domain $\Omega\subset\mathbb{R}^n$, $\forall 1\leq q,r\leq\infty$, and a natural number $m$. Suppose that a real number $\theta$ and a natural number $j$ are such that
\begin{align*}
\frac{1}{p}=\frac{j}{n}+(\frac{1}{r}-\frac{m}{n})\theta+\frac{1-\theta}{q},
\end{align*}
and
\begin{align*}
\frac{j}{m}\leq\theta\leq1,
\end{align*}
then there exists constant $C$ independent of $u$ such that
\begin{align*}
\|D^ju\|_{L^p(\Omega)}\leq C\|D^m u\|^\theta_{L^r(\Omega)}\|u\|^{1-\theta}_{L^q(\Omega)}.
\end{align*}
\end{lemma}

\section{Main result}

Before showing our main result, we give the definition of a weak solution for CNS and a strong solution for CPE. Recently, Bresch and Jabin \cite{br} consider different compactness method from Lions or Feireisl which can be applied to anisotropical stress tensor.\footnote{Let us  emphasize the result of Bresch and Jabin is valid only for small coefficients of viscosities.} They obtain the global existence of weak solutions for non-monotone pressure. Let us recall their definitions here.

\subsection{Dissipative weak solutions of CNS}

\begin{definition}\label{def1}
We say that $[\rho_\epsilon,\mathbf u_\epsilon]$ $\mathbf u_\epsilon=(\mathbf v_\epsilon,w_\epsilon)$ is a finite energy weak solution to the system of \eqref{aa}, supplemented with initial data (\ref{pa}) if $\rho_\epsilon=\rho_\epsilon(x,t)$ and
\begin{align}
&\mathbf u_\epsilon\in L^2(0,T;H^1(\Omega)),\hspace{3pt} \rho|\mathbf u_\epsilon|^2\in L^\infty(0,T; L^1(\Omega)),\nonumber\\
&\rho_\epsilon\in L^\infty(0,T;L^\gamma(\Omega))\cap C([0,T],L^1(\Omega)),
\end{align}

\noindent
$\bullet$ the continuity equation
\begin{align}
[\int_\Omega\rho_\epsilon\psi dxdz]^{t=\tau}_{t=0}=\int^\tau_0\int_{\Omega}\rho_\epsilon\partial_t\psi+\rho_\epsilon\mathbf{v}_\epsilon\cdot\nabla_x\psi+\rho_\epsilon w_\epsilon\partial_z\psi dxdzdt,
\end{align}
holds for all $\psi\in C^\infty_c([0,T)\times\Omega)$;

\noindent
$\bullet$
the momentum equation
\begin{align}
[\int_\Omega\rho_\epsilon\mathbf v_\epsilon\mathbf{\varphi_H} dxdz]^{t=\tau}_{t=0}-\int^\tau_0\int_{\Omega}\rho_\epsilon\mathbf{v}_\epsilon\partial_t\mathbf{\varphi_H}dxdzdt
-\int^\tau_0\int_{\Omega}\rho_\epsilon\mathbf{u}_\epsilon\mathbf{v}_\epsilon\cdot\nabla\mathbf{\varphi_H}dxdzdt\nonumber\\
+\int^\tau_0\int_\Omega\nabla\mathbf v_\epsilon: \nabla\mathbf{\varphi_H}dxdzdt-\int^\tau_0\int_\Omega p(\rho_\epsilon)\text{div}_x\mathbf{\varphi_H} dxdzdt=0,\label{3.3}
\end{align}
and
\begin{align}
\epsilon^2[\int_\Omega\rho_\epsilon w_\epsilon\varphi_3 dxdz]^{t=\tau}_{t=0}-\epsilon^2\int^\tau_0\int_{\Omega}\rho_\epsilon w_\epsilon\partial_t\varphi_3dxdzdt-\epsilon^2\int^\tau_0\int_{\Omega} \rho_\epsilon\mathbf{u}_\epsilon w_\epsilon\cdot\nabla\varphi_3dxdzdt\nonumber\\
+\epsilon^2\int^\tau_0\int_\Omega\nabla w_\epsilon\cdot\nabla\varphi_3dxdzdt-\int^\tau_0\int_\Omega p(\rho_\epsilon)\partial_z\varphi_3 dxdzdt=0,\label{3.4}
\end{align}
holds for all $\mathbf{\varphi_H},\varphi_3\in C^\infty_c([0,T)\times\Omega)$. Combining $(\ref{3.3})-(\ref{3.4})$, we obtain
\begin{align}
[\int_\Omega&\rho_\epsilon\mathbf v_\epsilon\mathbf{\varphi_H} dxdz
+\epsilon^2\int_\Omega\rho_\epsilon w_\epsilon\varphi_3 dxdz]^{t=\tau}_{t=0}\nonumber\\
&-\int^\tau_0\int_{\Omega}\rho_\epsilon\mathbf{v}_\epsilon\partial_t\mathbf{\varphi_H}dxdzdt-\epsilon^2\int^\tau_0\int_{\Omega}\rho_\epsilon w_\epsilon\partial_t\varphi_3dxdzdt\nonumber\\
&-\int^\tau_0\int_{\Omega}\rho_\epsilon\mathbf{v}_\epsilon\otimes\mathbf{v}_\epsilon:\nabla_x\mathbf{\varphi_H}dxdzdt
-\int^\tau_0\int_{\Omega}\rho_\epsilon\mathbf{v}_\epsilon w_\epsilon\cdot\partial_z\varphi_Hdxdzdt\nonumber\\
&-\epsilon^2\int^\tau_0\int_{\Omega} \rho_\epsilon\mathbf{v}_\epsilon w_\epsilon\cdot\nabla_x\varphi_3dxdzdt
-\epsilon^2\int^\tau_0\int_{\Omega} \rho_\epsilon w^2_\epsilon\partial_z\varphi_3dxdzdt\nonumber\\
&+\int^\tau_0\int_\Omega\nabla\mathbf v_\epsilon:\nabla\mathbf{\varphi_H}dxdzdt
+\epsilon^2\int^\tau_0\int_\Omega\nabla w_\epsilon\cdot\nabla\varphi_3dxdzdt-\int^\tau_0\int_\Omega p(\rho_\epsilon)\text{div}\mathbf{\varphi} dxdzdt=0,
\end{align}
where $\varphi=(\mathbf{\varphi_H},\varphi_3)\in C^\infty_c([0,T)\times\Omega)$ and $\text{div}\mathbf{\varphi}=\text{div}_x\mathbf{\varphi_H}+\partial_z\varphi_3$,

\noindent
$\bullet$
the energy inequality
\begin{align}
[\int_{\Omega}\frac{1}{2}\rho_\epsilon|\mathbf{v}_\epsilon|^2+\frac{\epsilon^2}{2}\rho_\epsilon |w_\epsilon|^2
+P(\rho_\epsilon))dxdz]|^{t=\tau}_{t=0}
+\int^\tau_0\int_\Omega(|\nabla\mathbf v_\epsilon|^2+\epsilon^2|\nabla w_\epsilon|^2)dxdzdt\leq 0,
\end{align}
holds for a.a $\tau\in(0,T)$, where $P(\rho)=\rho\int^\rho_1\frac{p(z)}{z^2}dz$.
\end{definition}
\subsection{Strong solution of CPE}

We say that $(r,\mathbf U)$, $\mathbf U=(\mathbf V,W)$ is a strong solution to the CPE system (\ref{aaa}) in $(0,T)\times\Omega$, if
\begin{align*}
&r^\frac{1}{2}\in L^\infty(0,T;H^2(\Omega)),\hspace{3pt}\partial_tr^\frac{1}{2}\in L^\infty(0,T;H^1(\Omega)),\hspace{3pt}r>0\hspace{3pt}\text{for all}\hspace{3pt}(t,x),\\
&\mathbf V\in L^\infty(0,T;H^3(\Omega))\cap L^2(0,T;H^4(\Omega)),\hspace{3pt} \partial_t\mathbf V\in L^2(0,T; H^2(\Omega)),
\end{align*}
with initial data $r^\frac{1}{2}_0\in H^2(\Omega)$, $r_0>0$ and $\mathbf V_0\in H^3(\Omega)$. Liu and Titi \cite{liu1} has proved the local existence of strong solution to CPE system (\ref{aaa}).

\begin{remark}
As the density is independent of $z$, we can obtain the following information of vertical velocity for the weak solution of CNS :
\begin{align}
\rho w(x,z,t)=-\rm{div}_x(\rho\widetilde{\mathbf v})+z\rm{div}_x(\rho\overline{\mathbf v}), \hspace{4pt}
\text{in the sense of} \hspace{4pt}H^{-1}(\Omega),
\label{b1}
\end{align}
where
\begin{align*}
\widetilde{\mathbf v}(x,z,t)=\int^z_0v(x,s,t)ds,\hspace{5pt}\overline{v}(x,t)=\int^1_0v(x,z,t)dz.
\end{align*}
Similarly, we can obtain the same equation for the strong solution of CPE in the classical sense. There is no information about $w$, so we need to derive its information. We should emphasize that \eqref{b1} is the key step to obtain the existence of weak solution for CPE in \cite{liu2,w}, which is inspired by incompressible case.
\end{remark}

\subsection{Relative entropy inequality}

Motivated by \cite{e2,e3}, for any finite energy weak solution $(\rho,\mathbf u)$, where $\mathbf u=(\mathbf v,w)$, to the CNS system, we introduce the relative energy functional
\begin{align}
\mathcal{E}(\rho,\mathbf{u}|r, \mathbf{U})&=\int_{\Omega}[\frac{1}{2}\rho|\mathbf v-\mathbf V|^2+\frac{\epsilon^2}{2}\rho|w-W|^2
+P(\rho)-P'(r)(\rho-r)-P(r)]dxdz\nonumber\\
&=\int_\Omega(\frac{1}{2}\rho|\mathbf v|^2+\frac{\epsilon^2}{2}\rho|w|^2+P(\rho))dxdz-\int_\Omega(\rho\mathbf v\cdot\mathbf V+\epsilon^2\rho wW)dxdz\nonumber\\
&\hspace{5pt}+\int_\Omega[\rho\frac{|\mathbf V|^2}{2}+\frac{\epsilon^2}{2}\rho|W|^2-\rho P'(r)]dxdz+\int_\Omega p(r)dxdz\nonumber\\
&=\sum^4_{i=1}I_i,\label{a1}
\end{align}
where $r>0$, $\mathbf U=(\mathbf V,W)$ are smooth ``test'' functions, $r$, $\mathbf U$ compactly supported in $\Omega$. Here we have used $rP'(r)-P(r)=p(r)$.

\begin{lemma}\label{relativeentropy}
Let $(\rho,\mathbf{v}, w)$ be a dissipative weak solution introduced in Definition  \ref{def1}. Then $(\rho,\mathbf{v}, w)$
satisfy the relative entropy inequality

\begin{align}
\mathcal{E}&(\rho,\mathbf{u}|r,\mathbf U)|^{t=\tau}_{t=0}+\int^\tau_0\int_\Omega\big{(}\nabla\mathbf v:(\nabla\mathbf v-\nabla\mathbf V)+\epsilon^2|\nabla w|^2\big{)}dxdzdt\nonumber\\
&\leq\int^\tau_0\int_{\Omega}\rho(\partial_t\mathbf V+\mathbf v\cdot\nabla_x\mathbf V+w\partial_z\mathbf V)(\mathbf V-\mathbf v)dxdt\nonumber\\
&\hspace{15pt}+\epsilon^2\int^\tau_0\int_{\Omega}\rho (\partial_tW+\mathbf v\cdot\nabla_xW+w\partial_zW)(W-w)dxdzdt
+\epsilon^2\int^\tau_0\int_\Omega\nabla w\cdot\nabla Wdxdzdt\nonumber\\
&\hspace{15pt}-\int^\tau_0\int_{\Omega}P''(r)((\rho-r)\partial_tr+\rho\mathbf v\cdot\nabla_xr)dxdzdt
-\int^\tau_0\int_{\Omega}p(r){\rm div}_x\mathbf Vdxdzdt.
\end{align}
\end{lemma}
{\bf Proof:}
From the weak formulation and energy inequality (3.3)-(3.6), we deduce
\begin{align}
&I_1|^{t=\tau}_{t=0}+\int^\tau_0\int_\Omega(|\nabla\mathbf v|^2+\epsilon^2|\nabla w|^2)dxdzdt\leq0,\\
&I_2|^{t=\tau}_{t=0}=-\int^\tau_0\int_\Omega\rho\mathbf v\partial_t\mathbf V+\rho\mathbf v\otimes\mathbf v:\nabla_x\mathbf v
+\rho\mathbf vw\cdot\partial_z\mathbf Vdxdzdt\nonumber\\
&\hspace{40pt}-\int^\tau_0\int_\Omega\epsilon^2\rho w\partial_tW+\epsilon^2\rho w(\mathbf v\cdot\nabla_x)W
+\epsilon^2\rho w^2\partial_zW+p(\rho)\text{div}_x\mathbf Vdxdzdt\nonumber\\
&\hspace{40pt}+\int^\tau_0\int_\Omega\nabla\mathbf v:\nabla\mathbf V+\epsilon^2\nabla w\cdot\nabla Wdxdzdt,\\
&I_3|^{t=\tau}_{t=0}=\int^\tau_0\int_\Omega\rho\partial_t\frac{|\mathbf V|^2}{2}+\rho\mathbf v\cdot\nabla_x\frac{|\mathbf V|^2}{2}+\rho w\partial_z\frac{|\mathbf V|^2}{2}dxdzdt\nonumber\\
&\hspace{40pt}+\epsilon^2\int^\tau_0\int_\Omega\rho\partial_t\frac{|W|^2}{2}+\rho\mathbf v\cdot\nabla_x\frac{|W|^2}{2}+\rho w\partial_z\frac{|W|^2}{2}dxdzdt\nonumber\\
&\hspace{40pt}-\int^\tau_0\int_\Omega\rho\partial_tP'(r)+\rho\mathbf v\cdot\nabla_xP'(r)+\rho w\partial_zP'(r)dxdzdt\nonumber\\
&\hspace{20pt}=\int^\tau_0\int_\Omega\rho\mathbf V\partial_t\mathbf V+\rho\mathbf v(\mathbf V\cdot\nabla_x)\mathbf V+\rho w\mathbf V\partial_z\mathbf Vdxdzdt\nonumber\\
&\hspace{30pt}+\epsilon^2\int^\tau_0\int_\Omega(\rho W\partial_tW+\rho W\mathbf v\cdot\nabla_xW+\rho wW\partial_zW)dxdzdt\nonumber\\
&\hspace{30pt}-\int^\tau_0\int_\Omega\rho P''(r)\partial_tr+P''(r)\rho\mathbf v\cdot\nabla_xr dxdzdt,\\
&I_4|^{t=\tau}_{t=0}=\int^\tau_0\int_\Omega\partial_tp(r)dxdzdt.
\end{align}

Summing (3.10)-(3.13) together, we obtain Lemma 3.1.

\begin{remark}
In the context of continuum mechanics, the concept of relative entropy has many successful applications, providing weak-strong uniqueness results, as well as justification of singular limits, for both incompressible and compressible fluid models. However, it is important to point out that there are distinguished differences for relative entropy inequality between incompressible and compressible cases. We suggest readers refer to \cite{e2,e3} for more details.
\end{remark}

\begin{remark}
Compared with the previous results \cite{gao}, there are some delicate differences in the process of using relative energy inequality. We should emphasize that we obtain the weak-strong uniqueness that is from weak solutions of CPE to strong solutions of CPE in \cite{gao}. Here, our convergence is between two different systems and is from 3D to 2.5D. The role of weak solutions is played by the solutions of CNS, and the strong solutions is played by those of CPE. It means that we should deal with the convergence of the vertical velocity of CNS and the absence of the information on the vertical velocity in CPE. Due to the special structure of CPE, the relative entropy inequality (3.8) is constructed differently from the ones in \cite{gao,ma}. Moreover, the adiabatic coefficient ($\gamma>4$) in the present work Theorem 3.1 which satisfies the assumptions of Bresch and Jabin's result \cite{br}($\gamma>\frac{3}{2}(\frac{4}{3}+\frac{\sqrt{10}}{3})\simeq3.5$), it improves our previous work \cite{gao}($\gamma>6$).
\end{remark}

Based on the relative entropy inequality, we can obtain the following lemma from \cite{e2}
\begin{lemma}
Let $0<a<b<\infty$. Then there exists $c=c(a,b)>0$ such that for all $\rho\in[0,\infty)$ and $r\in[a,b]$ there holds
\begin{equation} \label{pres}
P(\rho)-P'(r)(\rho-r)-P(r)\geq\left\{
\begin{array}{llll} C|\rho-r|^2,\hspace{5pt}\text{when} \hspace{3pt} \frac{r}{2}<\rho<r, \nonumber\\
C(1+\rho^\gamma),\hspace{5pt}\text{otherwise},
\end{array}\right.
\end{equation}
where $C=C(a,b)$.
\end{lemma}
Moreover, from \cite{e2}, we learn that
\begin{align}
&\mathcal{E}(\rho,\mathbf{u}|r,\mathbf U)(t)\in L^\infty(0,T),\hspace{3pt}
\int_\Omega\chi_{\rho\geq r}\rho^{\gamma}dxdz\leq C\mathcal{E}(\rho,\mathbf{u}|r,\mathbf U)(t),\nonumber\\
&\int_\Omega\chi_{\rho\leq \frac{r}{2}}1dxdz\leq C\mathcal{E}(\rho,\mathbf{u}|r,\mathbf U)(t),\hspace{3pt}
\int_\Omega\chi_{\frac{r}{2}<\rho<r}(\rho-r)^2dxdz\leq C\mathcal{E}(\rho,\mathbf{u}|r,\mathbf U)(t).\label{a3}
\end{align}
For a rigorous proof of Lemma 3.2 and (\ref{a3}), the reader is referred to \cite{e2}.

\subsection{Main result}
Now, we are ready to state our main result.
\begin{theorem}
Let $\gamma>4$, $T_{max}>0$ be the life time of strong solution to CPE system (\ref{aaa}) corresponding to initial data $[r_0,\mathbf V_0]$. Let $(\rho_\epsilon,\mathbf u_\epsilon)$, $\mathbf u_\epsilon=(\mathbf v_\epsilon,w_\epsilon)$ be a sequence of dissipative weak solutions to the CNS system $(\ref{aa})$ from the initial data $(\rho_{0,\epsilon},\mathbf u_{0,\epsilon})$. Suppose that
\begin{align*}
\mathcal{E}(\rho_{0,\epsilon},\mathbf{u}_{0,\epsilon}|r_0,\mathbf U_0)\rightarrow0,
\end{align*}
where $\mathbf U_0=(\mathbf V_0,W_0)$, then
\begin{align*}
ess\sup_{t\in(0,T_{max})}\mathcal{E}(\rho_{\epsilon},\mathbf{u}_{\epsilon}|r,\mathbf U)\rightarrow0,
\end{align*}
where $\mathbf U=(\mathbf V,W)$ and the couple $(r,\mathbf U)$ satisfy the CPE system (\ref{aaa}) on the time interval $[0,T_{max})$.

\end{theorem}

\begin{remark}
Recently, Bresch and Burtea \cite{b3} proved existence of weak solutions for anisotropic compressible Stokes system.
\end{remark}


Section 4 is devoted to the proof of the above theorem.

\section{Convergence}
In this section, we will prove the Theorem 3.1. First, we will explain the idea of the proof.
\subsection{Main idea of Proof}
The proof of Theorem 3.1 relies on the relative energy inequality by considering the strong solution $(r,\mathbf U)$, where $\mathbf U=(\mathbf V, W)$, as test function in the relative energy inequality \eqref{a1}. Firstly, let us recall the relative energy inequality
\begin{align}
\mathcal{E}&(\rho_{\epsilon},\mathbf{u}_{\epsilon}|r,\mathbf U)|^{t=\tau}_{t=0}+\int^\tau_0\int_\Omega\big{(}\nabla\mathbf v_{\epsilon}\cdot(\nabla\mathbf v_{\epsilon}-\nabla\mathbf V)+\epsilon^2\nabla w_{\epsilon}(\nabla w_{\epsilon}-\nabla W)\big{)}dxdzdt\nonumber\\
&\leq\int^\tau_0\int_{\Omega}\rho_{\epsilon}(\partial_t\mathbf V+\mathbf v_{\epsilon}\nabla_x\mathbf V+w_\epsilon\partial_z\mathbf V)(\mathbf V-\mathbf v_{\epsilon})dxdt\nonumber\\
&+\epsilon^2\int^\tau_0\int_{\Omega}\rho_{\epsilon} (\partial_tW+\mathbf v_{\epsilon}\nabla_xW+w_{\epsilon}\partial_zW)(W-w_{\epsilon})dxdzdt
+\epsilon^2\int^\tau_0\int_\Omega\nabla w_{\epsilon}\cdot\nabla Wdxdzdt\nonumber\\
&\hspace{15pt}-\int^\tau_0\int_{\Omega}P''(r)((\rho_{\epsilon}-r)\partial_tr+\rho_{\epsilon}\mathbf v_{\epsilon}\nabla_xr)dxdzdt
-\int^\tau_0\int_{\Omega}p(r){\rm div}_x\mathbf Vdxdzdt.\label{4a}
\end{align}
The goal now is to find a lower bound of (\ref{4a}) in the following form

\begin{align*}
\mathcal{E}(\rho_{\epsilon},\mathbf{u}_\epsilon|r,\mathbf U)(t)
+C\int^t_0\|\nabla\mathbf v_\epsilon-\nabla\mathbf V\|^2_{L^{2}}dt+\epsilon^2\int^t_0\|\nabla\mathbf w_\epsilon\|^2_{L^{2}}dt
\end{align*}
and an upper bound of (\ref{4a}) in the form
\begin{align*}
C(\delta)\int^t_0h(t)\mathcal{E}(\rho_{\epsilon},\mathbf{u}_\epsilon|r,\mathbf U)dt
+\delta\int^t_0\|\nabla\mathbf v_\epsilon-\nabla\mathbf V\|^2_{W^{1,2}}dt+o(\epsilon^2),
\end{align*}
with any $\delta>0$, where $C$ is independent of $\delta$ and $\epsilon$, $h\in L^1(0,T)$ and $o(\epsilon^2)\rightarrow0$ when $\epsilon\rightarrow0$.

Using the above bounds, we can deduce
\begin{align*}
\mathcal{E}(\rho_{\epsilon},\mathbf u_\epsilon|r,\mathbf U)(\tau)\leq C\int^\tau_0h(t)\mathcal{E}(\rho_{\epsilon},\mathbf u_{\epsilon}|r,\mathbf U)(t)dt+o(\epsilon^2),
\end{align*}
that implies the claim  by using the Gronwall inequality. In the rest of this section, we will follow this way.

\subsection{Step 1}
We write
\begin{align*}
\int_\Omega\rho_{\epsilon}&\mathbf v_{\epsilon}(\mathbf V-\mathbf v_{\epsilon})\cdot\nabla_x\mathbf Vdxdz=\\
&\int_\Omega\rho_{\epsilon}(\mathbf v_{\epsilon}-\mathbf V)(\mathbf V-\mathbf v_{\epsilon})\cdot\nabla_x\mathbf Vdxdz
+\int_\Omega\rho_{\epsilon}\mathbf V(\mathbf V-\mathbf v_{\epsilon})\cdot\nabla_x\mathbf Vdxdz.
\end{align*}

As $[r,\mathbf V, W]$ is a strong solution, it is obvious to obtain that
\begin{align}
\int_\Omega\rho_{\epsilon}(\mathbf v_{\epsilon}-\mathbf V)(\mathbf V-\mathbf v_{\epsilon})\cdot\nabla_x\mathbf Vdxdz
\leq C\mathcal{E}(\rho_{\epsilon},\mathbf u_{\epsilon}|r,\mathbf U).
\end{align}

Moreover, the momentum equation reads as
\begin{align*}
(r\mathbf V)_t+\text{div}_x(r\mathbf V\otimes\mathbf V)+\partial_z(r\mathbf VW)+\nabla_xp(r)=\Delta\mathbf V=\Delta_x\mathbf V+\partial_{zz}\mathbf V,
\end{align*}
which implies that
\begin{align*}
\mathbf V_t+\mathbf V\cdot\nabla_x\mathbf V+W\partial_z\mathbf V=-\frac{1}{r}\nabla_xp(r)+\frac{1}{r}\Delta_x\mathbf V
+\frac{1}{r}\partial_{zz}\mathbf V.
\end{align*}

So we rewrite the preceding two items on the right side of (\ref{4a}) as
\begin{align*}
\int_\Omega\rho_{\epsilon}[&\partial_t\mathbf V+\mathbf V\nabla_x\mathbf V+W\partial_z\mathbf V+(\mathbf v_{\epsilon}-\mathbf V)\nabla_x\mathbf V+(w_{\epsilon}-W)\partial_z\mathbf V](\mathbf V-\mathbf v_{\epsilon})dxdz\\
&=\int_\Omega\frac{\rho_{\epsilon}}{r}(\mathbf V-\mathbf v_{\epsilon})(\Delta_x\mathbf V+\partial_{zz}\mathbf V-\nabla_xp(r))dxdz\\
&\hspace{5pt}+\int_\Omega\rho_{\epsilon}(w_{\epsilon}-W)(\mathbf V-\mathbf v_{\epsilon})\cdot\partial_z\mathbf Vdxdz
-\int_\Omega\rho_{\epsilon}(\mathbf V-\mathbf v_{\epsilon})^2\nabla_x\mathbf V,
\end{align*}
and
\begin{align}
\epsilon^2\int^\tau_0&\int_{\Omega}\rho_{\epsilon} (\partial_tW+\mathbf v_{\epsilon}\nabla_xW+w_{\epsilon}\partial_zW)(W-w_{\epsilon})dxdzdt\nonumber\\
&\leq\int^\tau_0\mathcal{E}(\rho_{\epsilon},\mathbf{u}_{\epsilon}|r,\mathbf U)
+\epsilon^4\int^\tau_0\int_{\Omega}\rho_{\epsilon}(\partial_tW+\mathbf v_{\epsilon}\nabla_xW+w_{\epsilon}\partial_zW)^2dxdzdt\nonumber\\
&=\int^\tau_0\mathcal{E}(\rho_{\epsilon},\mathbf{u}_{\epsilon}|r,\mathbf U)
+\epsilon^4\int^\tau_0\int_{\Omega}\rho_{\epsilon}(\partial_tW+\mathbf V\nabla_xW+W\partial_zW)^2dxdzdt\nonumber\\
&\hspace{10pt}+\epsilon^4\int^\tau_0\int_{\Omega}\rho_{\epsilon}((\mathbf v_{\epsilon}-\mathbf V)\nabla_xW+(w_{\epsilon}-W)\partial_zW)^2dxdzdt.\label{4.3}
\end{align}
Noticing Lemma 3.2, we have
\begin{align}
\int^\tau_0\int_{\Omega}&\rho_{\epsilon}(\partial_tW+\mathbf V\nabla_xW+W\partial_zW)^2dxdzdt\nonumber\\
&=\int^\tau_0\int_{\Omega}\chi_{\rho_{\epsilon}< \frac{r}{2}}\rho_{\epsilon}(\partial_tW+\mathbf V\nabla_xW+W\partial_zW)^2dxdz\nonumber\\
&\hspace{5pt}+\int^\tau_0\int_{\Omega}\chi_{\frac{r}{2}\leq\rho_{\epsilon}\leq r}\rho_{\epsilon}(\partial_tW+\mathbf V\nabla_xW+W\partial_zW)^2dxdz\nonumber\\
&\hspace{5pt}+\int^\tau_0\int_{\Omega}\chi_{\rho_{\epsilon}>r}\rho_{\epsilon}(\partial_tW+\mathbf V\nabla_xW+W\partial_zW)^2dxdzdt\nonumber\\
&\leq\int^\tau_0\int_{\Omega}\chi_{\rho_{\epsilon}< \frac{r}{2}}r(\partial_tW+\mathbf V\nabla_xW+W\partial_zW)^2dxdz\nonumber\\
&\hspace{5pt}+\int^\tau_0\int_{\Omega}\chi_{\rho_{\epsilon}>r}\rho_{\epsilon}(\partial_tW+\mathbf V\nabla_xW+W\partial_zW)^2dxdzdt\nonumber\\
&\hspace{5pt}+C\int^\tau_0\int_{\Omega}\chi_{\frac{r}{2}\leq\rho_{\epsilon}\leq r}(\rho_{\epsilon}-r)(\partial_tW+\mathbf V\nabla_xW+W\partial_zW)^2dxdz+C\nonumber\\
&\leq C\int^\tau_0\mathcal{E}(\rho_{\epsilon},\mathbf{u}_{\epsilon}|r,\mathbf U)+
C\int^\tau_0\int_{\Omega}\chi_{\frac{r}{2}\leq\rho_{\epsilon}\leq r}(\rho_{\epsilon}-r)^2dxdz+C\int^\tau_0\int_{\Omega}\chi_{\rho_{\epsilon}>r}\rho_{\epsilon}^\gamma dxdz+C\nonumber\\
&\leq C\int^\tau_0\mathcal{E}(\rho_{\epsilon},\mathbf{u}_{\epsilon}|r,\mathbf U)+C.\label{4.4}
\end{align}
It is easy to use the definition of relative entropy inequality and Cauchy inequality to obtain
\begin{align}
\epsilon^4\int^\tau_0\int_{\Omega}\rho_{\epsilon}((\mathbf v_{\epsilon}-\mathbf V)\nabla_xW+(w_{\epsilon}-W)\partial_zW)^2dxdzdt
\leq \epsilon^2\mathcal{E}(\rho_{\epsilon},\mathbf{u}_{\epsilon}|r,\mathbf U)+o(\epsilon^2).\label{444}
\end{align}

Putting (\ref{4.4})-(\ref{444}) into (\ref{4.3}) yields
\begin{align*}
\epsilon^2\int^\tau_0&\int_{\Omega}\rho_{\epsilon} (\partial_tW+\mathbf v_{\epsilon}\nabla_xW+w_{\epsilon}\partial_zW)(W-w_{\epsilon})dxdzdt
\leq C\int^\tau_0\mathcal{E}(\rho_{\epsilon},\mathbf{u}_{\epsilon}|r,\mathbf U)+o(\epsilon^2).
\end{align*}

Moreover, a simple application of Cauchy inequality leads to the following
\begin{align*}
\epsilon^2\int^\tau_0\int_\Omega\nabla w_{\epsilon}\cdot\nabla Wdxdzdt\leq\frac{\epsilon^2}{2}\int^\tau_0\int_\Omega|\nabla w_{\epsilon}|^2dxdzdt+o(\epsilon^2).
\end{align*}

Thus, we obtain that
\begin{align*}
\mathcal{E}&(\rho_{\epsilon},\mathbf{u}_{\epsilon}|r,\mathbf U)|^{t=\tau}_{t=0}+\int^\tau_0\int_\Omega\big{(}\nabla\mathbf v_{\epsilon}\cdot(\nabla\mathbf v_{\epsilon}-\nabla\mathbf V)+\frac{\epsilon^2}{2}|\nabla w_{\epsilon}|^2\big{)}dxdzdt\nonumber\\
&\leq C\int^\tau_0\mathcal{E}(\rho_{\epsilon},\mathbf{u}_{\epsilon}|r,\mathbf U)dt
-\int^\tau_0\int_{\Omega}P''(r)((\rho_{\epsilon}-r)\partial_tr+\rho_{\epsilon}\mathbf v_{\epsilon}\nabla_xr)dxdzdt\nonumber\\
&\hspace{8pt}+\int^\tau_0\int_\Omega\frac{\rho_{\epsilon}}{r}(\mathbf V-\mathbf v_{\epsilon})(\Delta_x\mathbf V+\partial_{zz}\mathbf V)dxdz
-\int^\tau_0\int_\Omega\frac{\rho_{\epsilon}}{r}(\mathbf V-\mathbf v_{\epsilon})\nabla_xp(r)dxdz\nonumber\\
&\hspace{8pt}+\int^\tau_0\int_\Omega\rho_{\epsilon}(w_{\epsilon}-W)(\mathbf V-\mathbf v_{\epsilon})\cdot\partial_z\mathbf Vdxdzdt
-\int^\tau_0\int_\Omega p(\rho_{\epsilon})\text{div}_x\mathbf Vdxdzdt+o(\epsilon^2).\nonumber\\
\end{align*}

\subsection{Step 2}
The major challenges of the analysis is to estimate the complicated nonlinear term $\int_\Omega\rho_{\epsilon}(w_{\epsilon}-W)(\mathbf V-\mathbf v_{\epsilon})\cdot\partial_z\mathbf Vdxdz$, we rewrite it as
\begin{align}
\int_\Omega\rho_{\epsilon}&(w_{\epsilon}-W)(\mathbf V-\mathbf v_{\epsilon})\cdot\partial_z\mathbf Vdxdz\nonumber\\
&=\int_\Omega\rho_{\epsilon} w_{\epsilon}(\mathbf V-\mathbf v_{\epsilon})\cdot\partial_z\mathbf Vdxdz-\int_\Omega\rho_{\epsilon} W(\mathbf V-\mathbf v_{\epsilon})\cdot\partial_z\mathbf Vdxdz.\label{b}
\end{align}

A similar heuristic argument from \cite{e2,kr} shows that the second term on the right side of (\ref{b}) will be split into three parts
\begin{align}
\int_\Omega&\rho_{\epsilon} W(\mathbf V-\mathbf v_{\epsilon})\cdot\partial_z\mathbf Vdxdz\nonumber\nonumber\\
&=\int_\Omega\chi_{\rho_{\epsilon}\leq \frac{r}{2}}\rho_{\epsilon} W(\mathbf V-\mathbf v_{\epsilon})\cdot\partial_z\mathbf Vdxdz
+\int_\Omega\chi_{\frac{r}{2}<\rho_{\epsilon}<r}\rho_{\epsilon} W(\mathbf V-\mathbf v_{\epsilon})\cdot\partial_z\mathbf Vdxdz\nonumber\\
&\hspace{8pt}+\int_\Omega\chi_{\rho_{\epsilon}\geq r}\rho_{\epsilon} W(\mathbf V-\mathbf v_{\epsilon})\cdot\partial_z\mathbf Vdxdz\nonumber\nonumber\\
&\leq \|\chi_{\rho_{\epsilon}\leq \frac{r}{2}}1\|_{L^2(\Omega)}\|r\|_{L^\infty}\|W\partial_z\mathbf V\|_{L^3}\|\mathbf V-\mathbf v_{\epsilon}\|_{L^6(\Omega)}
+\int_\Omega\chi_{\rho_{\epsilon}\geq r}\rho_{\epsilon}^{\frac{\gamma}{2}}W\partial_z\mathbf V\cdot(\mathbf V-\mathbf v_{\epsilon})dxdz\nonumber\\
&\hspace{8pt}+C\|\chi_{\frac{r}{2}<\rho_{\epsilon}< r}(\rho_{\epsilon}-r)\|_{L^2(\Omega)}\|W\partial_z\mathbf V\|_{L^3}\|\mathbf V-\mathbf v_{\epsilon}\|_{L^6(\Omega)}\nonumber\\
&\leq C\int_\Omega\chi_{\rho_{\epsilon}\leq \frac{r}{2}}1dxdz +C\int_\Omega\chi_{\frac{r}{2}<\rho_{\epsilon}<r}(\rho_{\epsilon}-r)^2dxdz\nonumber\\
&\hspace{10pt}+C\int_\Omega\chi_{\rho_{\epsilon}\geq r}\rho_{\epsilon}^\gamma dxdz
+\delta\|\mathbf V-\mathbf v_{\epsilon}\|^2_{L^6(\Omega)}\nonumber\\
&\leq C\mathcal{E}(\rho_{\epsilon},\mathbf u_{\epsilon}|r,\mathbf U)+\delta\|\nabla_x\mathbf V-\nabla_x\mathbf v_{\epsilon}\|^2_{L^2(\Omega)}
+\delta\|\partial_z\mathbf V-\partial_z\mathbf v_{\epsilon}\|^2_{L^2(\Omega)},\label{33c}
\end{align}
where in the last inequality, we have used Lemma 2.1.

We now turn to analyze the first term on the right hand of \eqref{b}, which is the crucial and difficult part in our proof. Taking \eqref{b1} into it, we have
\begin{align}
\int_\Omega\rho_{\epsilon} w_{\epsilon}&(\mathbf V-\mathbf v_{\epsilon})\cdot\partial_z\mathbf Vdxdz\nonumber\\
&=\int_\Omega[-\text{div}_x(\rho_{\epsilon}\widetilde{\mathbf v}_{\epsilon})+z\text{div}_x(\rho_{\epsilon}\overline{\mathbf v}_{\epsilon})]\partial_z\mathbf V
\cdot(\mathbf V-\mathbf v_{\epsilon})dxdz\nonumber\\
&=\int_\Omega(\rho_{\epsilon}\widetilde{\mathbf v}_{\epsilon}-z\rho_{\epsilon}\overline{\mathbf v}_{\epsilon})\partial_z\nabla_x\mathbf V\cdot(\mathbf V-\mathbf v_{\epsilon})dxdz\nonumber\\
&\hspace{10pt}+\int_\Omega(\rho_{\epsilon}\widetilde{\mathbf v}_{\epsilon}-z\rho_{\epsilon}\overline{\mathbf v}_{\epsilon})\partial_z\mathbf V\cdot(\nabla_x\mathbf V-\nabla_x\mathbf v_{\epsilon})dxdz.\label{c}
\end{align}

In the following, we will estimate the terms on the right hand side of \eqref{c}. We only need to consider the most complicated terms, the remaining terms can be completed by the similar way. Firstly, we deal with $\int_\Omega\rho_{\epsilon}\widetilde{\mathbf v}_{\epsilon}\partial_z\nabla_x\mathbf V\cdot(\mathbf V-\mathbf v_{\epsilon})dxdz$ as the follows,
\begin{align*}
\int_\Omega\rho_{\epsilon}&\widetilde{\mathbf v}_{\epsilon}\partial_z\nabla_x\mathbf V\cdot(\mathbf V-\mathbf v_{\epsilon})dxdz\\
&=\int_\Omega\rho_{\epsilon}(\widetilde{\mathbf v}_{\epsilon}-\widetilde{\mathbf V})\partial_z\nabla_x\mathbf V\cdot(\mathbf V-\mathbf v_{\epsilon})dxdz
+\int_\Omega\rho_{\epsilon}\widetilde{\mathbf V}\partial_z\nabla_x\mathbf V\cdot(\mathbf V-\mathbf v_{\epsilon})dxdz\\
&=J_1+J_2,
\end{align*}
where $\widetilde{\mathbf V}=\int^z_0\mathbf V(x,s,t)ds$.

Similar to the above analysis, we decompose the term $J_2$ into three parts
\begin{align*}
J_2&=\int_\Omega\rho_{\epsilon}\widetilde{\mathbf V}\partial_z\nabla_x\mathbf V\cdot(\mathbf V-\mathbf v_{\epsilon})dxdz\\
&=\int_\Omega\chi_{\rho_{\epsilon}\leq \frac{r}{2}}\rho_{\epsilon}\widetilde{\mathbf V}\partial_z\nabla_x\mathbf V\cdot(\mathbf V-\mathbf v_{\epsilon})dxdz
+\int_\Omega\chi_{\frac{r}{2}<\rho_{\epsilon}<r}\rho_{\epsilon}\widetilde{\mathbf V}\partial_z\nabla_x\mathbf V\cdot(\mathbf V-\mathbf v_{\epsilon})dxdz\\
&\hspace{5pt}+\int_\Omega\chi_{\rho_{\epsilon}\geq r}\rho_{\epsilon}\widetilde{\mathbf V}\partial_z\nabla_x\mathbf V\cdot(\mathbf V-\mathbf v_{\epsilon})dxdz\\
&\leq \|\chi_{\rho_{\epsilon}\leq \frac{r}{2}}1\|_{L^2(\Omega)}\|r\|_{L^\infty}\|\widetilde{\mathbf V}\partial_z\nabla_x\mathbf V\|_{L^3}\|\mathbf V-\mathbf v_{\epsilon}\|_{L^6(\Omega)}\\
&\hspace{10pt}+\|\chi_{\rho_{\epsilon}\geq r}\rho_{\epsilon}^{\frac{\gamma}{2}}\|_{L^2(\Omega)}\|\widetilde{\mathbf V}\partial_z\nabla_x\mathbf V\|_{L^3(\Omega)}\|\mathbf V-\mathbf v_{\epsilon}\|_{L^6(\Omega)}\\
&\hspace{8pt}+C\|\chi_{\frac{r}{2}<\rho_{\epsilon}<r}(\rho_{\epsilon}-r)\|_{L^2(\Omega)}\|\widetilde{\mathbf V}\partial_z\nabla_x\mathbf V\|_{L^3(\Omega)}\|\mathbf V-\mathbf v_{\epsilon}\|_{L^6(\Omega)}\\
&\leq C\mathcal{E}(\rho_{\epsilon},\mathbf u_{\epsilon}|r,\mathbf U)(t)+\delta\|\nabla_x\mathbf V-\nabla_x\mathbf v_{\epsilon}\|^2_{L^2(\Omega)}
+\delta\|\partial_z\mathbf V-\partial_z\mathbf v_{\epsilon}\|^2_{L^2(\Omega)}.
\end{align*}

On the other hand, by virtue of Cauchy inequality, it follows that
\begin{align}
J_1&=\int_\Omega\rho_{\epsilon}(\widetilde{\mathbf v}_{\epsilon}-\widetilde{\mathbf V})\partial_z\nabla_x\mathbf V\cdot(\mathbf V-\mathbf v_{\epsilon})dxdz\nonumber\\
&\leq\|\partial_z\nabla_x\mathbf V\|_{L^\infty}\int_\Omega\rho_{\epsilon}|\widetilde{\mathbf v}_{\epsilon}-\widetilde{\mathbf V}|^2dxdz+\int_\Omega\rho_{\epsilon}|\mathbf V-\mathbf v_{\epsilon}|^2dxdz\nonumber\\
&\leq C\int_\Omega\rho_{\epsilon}|\int^z_0(\mathbf v_{\epsilon}(s)-\mathbf V(s))ds|^2dxdz+\mathcal{E}(\rho_{\epsilon},\mathbf u_{\epsilon}|r,\mathbf U)\nonumber\\
&\leq C\int_\Omega\rho_{\epsilon}\big{(}\int^1_0|\mathbf V-\mathbf v_{\epsilon}|^2ds\big{)}dxdz+\mathcal{E}(\rho_{\epsilon},\mathbf u_{\epsilon}|r,\mathbf U)\nonumber\\
&\leq C\int^1_0\int_\Omega \rho_{\epsilon}|\mathbf V-\mathbf v_{\epsilon}|^2dxdzds+\mathcal{E}(\rho_{\epsilon},\mathbf u_{\epsilon}|r,\mathbf U)\nonumber\\
&\leq C\int_\Omega \rho_{\epsilon}|\mathbf V-\mathbf v_{\epsilon}|^2dxdz+\mathcal{E}(\rho_{\epsilon},\mathbf u_{\epsilon}|r,\mathbf U)\nonumber\\
&\leq C\mathcal{E}(\rho_{\epsilon},\mathbf u_{\epsilon}|r,\mathbf U).\label{3aaa}
\end{align}

Secondly, we will investigate another complicated nonlinear term $\int_\Omega\rho_{\epsilon}\widetilde{\mathbf v}_\epsilon\partial_z\mathbf V\cdot(\nabla_x\mathbf V-\nabla_x\mathbf v_{\epsilon})dxdz$. It is straightforward to show that
\begin{align}
\int_\Omega\rho_{\epsilon}\widetilde{\mathbf v}_\epsilon&\partial_z\mathbf V\cdot(\nabla_x\mathbf V-\nabla_x\mathbf v_{\epsilon})dxdz\nonumber\\
&=\int_\Omega\chi_{\rho_{\epsilon}< r}\rho_{\epsilon}\widetilde{\mathbf v}_\epsilon\partial_z\mathbf V\cdot(\nabla_x\mathbf V-\nabla_x\mathbf v_{\epsilon})dxdz
+\int_\Omega\chi_{\rho_{\epsilon}\geq r}\rho_{\epsilon}\widetilde{\mathbf v}_\epsilon\partial_z\mathbf V\cdot(\nabla_x\mathbf V-\nabla_x\mathbf v_{\epsilon})dxdz,\label{2a}
\end{align}
where the first term on the right side of \eqref{2a} is split into two parts as
\begin{align*}
\int_\Omega\chi_{\rho_{\epsilon}< r}&\rho_{\epsilon}\widetilde{\mathbf v}_\epsilon\partial_z\mathbf V\cdot(\nabla_x\mathbf V-\nabla_x\mathbf v_{\epsilon})dxdz\\
&=\int_{\Omega}\chi_{\rho_{\epsilon}<r}\rho_{\epsilon}(\widetilde{\mathbf v}_{\epsilon}-\widetilde{\mathbf V})\partial_z\mathbf V\cdot(\nabla_x\mathbf V-\nabla_x\mathbf v_{\epsilon})dxdz\\
&\hspace{20pt}+\int_{\Omega}\chi_{\rho_{\epsilon}<r}\rho_{\epsilon}\widetilde{\mathbf V}\partial_z\mathbf V\cdot(\nabla_x\mathbf V-\nabla_x\mathbf v_{\epsilon})dxdz\\
&=\int_{\Omega}\chi_{\rho_{\epsilon}<r}\rho_{\epsilon}(\widetilde{\mathbf v}_{\epsilon}-\widetilde{\mathbf V})\partial_z\mathbf V\cdot(\nabla_x\mathbf V-\nabla_x\mathbf v_{\epsilon})dxdz\\
&\hspace{20pt}+\int_{\Omega}\chi_{\frac{r}{2}<\rho_{\epsilon}<r}\rho_{\epsilon}\widetilde{\mathbf V}\partial_z\mathbf V\cdot(\nabla_x\mathbf V-\nabla_x\mathbf v_{\epsilon})dxdz\\
&\hspace{20pt}+\int_{\Omega}\chi_{\rho_{\epsilon}\leq \frac{r}{2}}\rho_{\epsilon}\widetilde{\mathbf V}\partial_z\mathbf V\cdot(\nabla_x\mathbf V-\nabla_x\mathbf v_{\epsilon})dxdz\\
&\leq \|\chi_{\rho_{\epsilon}<r}\rho_{\epsilon}^{\frac{1}{2}}\|_{L^\infty(\Omega)}\|\sqrt{\rho_{\epsilon}}(\widetilde{\mathbf v}_{\epsilon}-\widetilde{\mathbf V})\|_{L^2(\Omega)}\|\partial_z\mathbf V\|_{L^\infty(\Omega)}\|\nabla_x\mathbf V-\nabla_x\mathbf v_{\epsilon}\|_{L^2(\Omega)}\\
&\hspace{10pt}+\|\chi_{\frac{r}{2}<\rho_{\epsilon}< r}\rho_{\epsilon}\|_{L^2(\Omega)}\|\widetilde{\mathbf V}\partial_z\mathbf V\|_{L^\infty(\Omega)}\|\nabla_x\mathbf V-\nabla_x\mathbf v_{\epsilon}\|_{L^2(\Omega)}\\
&\hspace{10pt}+\|\chi_{\rho_{\epsilon}\leq \frac{r}{2}}1\|_{L^2(\Omega)}\|r\|_{L^\infty(\Omega)}
\|\widetilde{\mathbf V}\partial_z\mathbf V\|_{L^\infty(\Omega)}\|\nabla_x\mathbf V-\nabla_x\mathbf v_{\epsilon}\|_{L^2(\Omega)}\\
&\leq C\mathcal{E}(\rho_{\epsilon},\mathbf u_{\epsilon}|r,\mathbf U)(t)+\delta\|\nabla_x\mathbf V-\nabla_x\mathbf v_{\epsilon}\|^2_{L^2(\Omega)}.
\end{align*}

The decomposition of remainder of \eqref{2a} is identical to the above as:
\begin{align}
\int_\Omega&\chi_{\rho_{\epsilon}\geq r}\rho_{\epsilon}\widetilde{\mathbf v}_\epsilon\partial_z\mathbf V\cdot(\nabla_x\mathbf V-\nabla_x\mathbf v_{\epsilon})dxdz\nonumber\\
&=\int_\Omega\chi_{\rho_{\epsilon}\geq r}\rho_{\epsilon}(\widetilde{\mathbf v}_{\epsilon}-\widetilde{\mathbf V})\partial_z\mathbf V\cdot(\nabla_x\mathbf V-\nabla_x\mathbf v_{\epsilon})dxdz
+\int_\Omega\chi_{\rho_{\epsilon}\geq r}\rho_{\epsilon}\widetilde{\mathbf V}\partial_z\mathbf V\cdot(\nabla_x\mathbf V-\nabla_x\mathbf v_{\epsilon})dxdz\nonumber\\
&=K_1+K_2,\label{333}
\end{align}
where
\begin{align}
K_2&\leq \int_\Omega\chi_{\rho_{\epsilon}\geq r}\rho_{\epsilon}^{\frac{\gamma}{2}}\widetilde{\mathbf V}\partial_z\mathbf V\cdot(\nabla_x\mathbf V-\nabla_x\mathbf v_{\epsilon})dxdz\nonumber\\
&\leq \|\chi_{\rho_{\epsilon}\geq r}\rho_{\epsilon}^{\frac{\gamma}{2}}\|_{L^2(\Omega)}
\|\widetilde{\mathbf V}\partial_z\mathbf V\|_{L^\infty(\Omega)}
\|\nabla_x\mathbf V-\nabla_x\mathbf v_{\epsilon}\|_{L^2(\Omega)}\nonumber\\
&\leq C\|\chi_{\rho_{\epsilon}\geq r}\rho_{\epsilon}^{\frac{\gamma}{2}}\|^2_{L^2(\Omega)}
+\delta\|\nabla_x\mathbf V-\nabla_x\mathbf v_{\epsilon}\|^2_{L^2(\Omega)}\nonumber\\
&\leq C\mathcal{E}(\rho_{\epsilon},\mathbf u_{\epsilon}|r,\mathbf U)(t)+\delta\|\nabla_x\mathbf V-\nabla_x\mathbf v_{\epsilon}\|^2_{L^2(\Omega)}.
\end{align}

It remains to estimate $K_1$. Due to H\"{o}lder inequality, it follows that
\begin{align*}
K_1&\leq \|\chi_{\rho_{\epsilon}\geq r}\rho_{\epsilon}\|_{L^4(\Omega)}
\|\chi_{\rho_{\epsilon}\geq r}(\widetilde{\mathbf v}_{\epsilon}-\widetilde{\mathbf V})\|_{L^4(\Omega)}
\|\partial_z\mathbf V\|_{L^\infty(\Omega)}
\|\nabla_x\mathbf V-\nabla_x\mathbf v_{\epsilon}\|_{L^2(\Omega)}\nonumber\\
&\leq C\|\chi_{\rho_{\epsilon}\geq r}\rho_{\epsilon}\|^2_{L^4(\Omega)}
\|\chi_{\rho_{\epsilon}\geq r}(\widetilde{\mathbf v}_{\epsilon}-\widetilde{\mathbf V})\|^2_{L^4(\Omega)}
+\delta\|\nabla_x\mathbf V-\nabla_x\mathbf v_{\epsilon}\|^2_{L^2(\Omega)}\nonumber\\
&\leq  C\|\chi_{\rho_{\epsilon}\geq r}\rho_{\epsilon}\|^2_{L^4(\Omega)}
\|\chi_{\rho_{\epsilon}\geq r}(\widetilde{\mathbf v}_{\epsilon}-\widetilde{\mathbf V})\|_{L^3(\Omega)}
\|\chi_{\rho_{\epsilon}\geq r}(\nabla\widetilde{\mathbf v}_{\epsilon}-\nabla\widetilde{\mathbf V})\|_{L^2(\Omega)}
+\delta\|\nabla_x\mathbf v_{\epsilon}-\nabla_x\mathbf V\|^2_{L^2(\Omega)}\nonumber\\
&\leq C\|\chi_{\rho_{\epsilon}\geq r}\rho_{\epsilon}\|^4_{L^4(\Omega)}
\|\chi_{\rho_{\epsilon}\geq r}(\widetilde{\mathbf v}_{\epsilon}-\widetilde{\mathbf V})\|^2_{L^3(\Omega)}
+\delta\|\nabla_x\widetilde{\mathbf v}_{\epsilon}-\nabla_x\widetilde{\mathbf V}\|^2_{L^2(\Omega)}\nonumber\\
&\hspace{10pt}+\delta\|\partial_z\widetilde{\mathbf v}_{\epsilon}-\partial_z\widetilde{\mathbf V}\|^2_{L^2(\Omega)}
+\delta\|\nabla_x\mathbf v_{\epsilon}-\nabla_x\mathbf V\|^2_{L^2(\Omega)}\nonumber\\
&\leq\|\chi_{\rho_{\epsilon}\geq r}\rho_{\epsilon}\|^4_{L^4(\Omega)}
\|\chi_{\rho_{\epsilon}\geq r}(\widetilde{\mathbf v}_{\epsilon}-\widetilde{\mathbf V})\|_{L^2(\Omega)}\|\chi_{\rho_{\epsilon}\geq r}(\widetilde{\mathbf v}_{\epsilon}-\widetilde{\mathbf V})\|_{H^1(\Omega)}
+\delta\|\nabla_x\widetilde{\mathbf v}_{\epsilon}-\nabla_x\widetilde{\mathbf V}\|^2_{L^2(\Omega)}\nonumber\\
&\hspace{10pt}+\delta\|\partial_z\widetilde{\mathbf v}_{\epsilon}-\partial_z\widetilde{\mathbf V}\|^2_{L^2(\Omega)}
+\delta\|\nabla_x\mathbf v_{\epsilon}-\nabla_x\mathbf V\|^2_{L^2(\Omega)}\nonumber\\
&\leq\|\chi_{\rho_{\epsilon}\geq r}\rho_{\epsilon}\|^8_{L^4(\Omega)}
\|\chi_{\rho_{\epsilon}\geq r}(\widetilde{\mathbf v}_{\epsilon}-\widetilde{\mathbf V})\|^2_{L^2(\Omega)}+\delta\|\chi_{\rho_{\epsilon}\geq r}(\widetilde{\mathbf v}_{\epsilon}-\widetilde{\mathbf V})\|^2_{L^2(\Omega)}
+\delta\|\nabla_x\widetilde{\mathbf v}_{\epsilon}-\nabla_x\widetilde{\mathbf V}\|^2_{L^2(\Omega)}\nonumber\\
&\hspace{10pt}+\delta\|\partial_z\widetilde{\mathbf v}_{\epsilon}-\partial_z\widetilde{\mathbf V}\|^2_{L^2(\Omega)}
+\delta\|\nabla_x\mathbf v_{\epsilon}-\nabla_x\mathbf V\|^2_{L^2(\Omega)}
\end{align*}
where we have used the Lemma 2.2
\begin{align*}
\|f\|_{L^4}\leq\|\nabla f\|^{\frac{1}{2}}_{L^2}\|f\|^{\frac{1}{2}}_{L^3}\hspace{3pt}{\rm and}\hspace{3pt}\|f\|_{L^3}\leq\|f\|^{\frac{1}{2}}_{L^2}\|f\|^{\frac{1}{2}}_{H^1}.
\end{align*}
Recalling \eqref{a3} and \eqref{3aaa}, we have
\begin{align*}
\|\chi_{\rho_{\epsilon}\geq r}\rho_{\epsilon}\|^8_{L^4(\Omega)}
=(\int_{\rho_{\epsilon}\geq r}\rho_{\epsilon}^4 dxdz)^2\leq C(\int_{\Omega}\rho_{\epsilon}^\gamma dxdz)^{\frac{8}{\gamma}}
\leq\mathcal{E}(\rho_{\epsilon},\mathbf u_{\epsilon}|r,\mathbf U)^{\frac{8}{\gamma}}(t),
\end{align*}
and
\begin{align*}
\|\chi_{\rho_{\epsilon}\geq r}(\widetilde{\mathbf v}_{\epsilon}-\widetilde{\mathbf V})\|^2_{L^2(\Omega)}
&=\int_{\rho_{\epsilon}\geq r}|\widetilde{\mathbf v}_{\epsilon}-\widetilde{\mathbf V}|^2dxdz
=\int_{\rho_{\epsilon}\geq r}\frac{1}{\rho_{\epsilon}}\rho_{\epsilon}|\widetilde{\mathbf v}_{\epsilon}-\widetilde{\mathbf V}|^2dxdz\\
&\leq \frac{1}{\|r\|_{\infty(\Omega)}}\mathcal{E}(\rho_{\epsilon},\mathbf u_{\epsilon}|r,\mathbf U)(t).
\end{align*}

An argument similar to the one used in \eqref{3aaa} yields
\begin{align*}
\|\nabla_x\widetilde{\mathbf v}_{\epsilon}-\nabla_x\widetilde{\mathbf V}\|^2_{L^2(\Omega)}
\leq \|\nabla_x\mathbf v_{\epsilon}-\nabla_x\mathbf V\|^2_{L^2(\Omega)},
\hspace{5pt}\|\partial_z\widetilde{\mathbf v}_{\epsilon}-\partial_z\widetilde{\mathbf V}\|^2_{L^2(\Omega)}
\leq\|\partial_z\mathbf v_{\epsilon}-\partial_z\mathbf V\|^2_{L^2(\Omega)}.
\end{align*}


Combining the above estimates, we arrive at the conclusion that
\begin{align*}
\int^\tau_0K_1dt\leq C\int^\tau_0h(t)\mathcal{E}(\rho_{\epsilon},\mathbf u_{\epsilon}|r,\mathbf U)(t)dt
+\delta\int^\tau_0\|\nabla_x\mathbf v_{\epsilon}-\nabla_x\mathbf V\|^2_{L^2(\Omega)}+
\|\partial_z\mathbf v_{\epsilon}-\partial_z\mathbf V\|^2_{L^2(\Omega)}dt,
\end{align*}
where $h(t)\in L^1(0,T)$.

The estimate of remainder in \eqref{c} can be completed by the analogous method. Therefore, we can summarize what we have proved as the following
\begin{align*}
\mathcal{E}&(\rho_{\epsilon},\mathbf{u}_{\epsilon}|r,\mathbf U)|^{t=\tau}_{t=0}+\int^\tau_0\int_\Omega\big{(}\nabla\mathbf v_{\epsilon}\cdot(\nabla\mathbf v_{\epsilon}-\nabla\mathbf V)+\epsilon^2|\nabla w_{\epsilon}|^2\big{)}dxdzdt\nonumber\\
&\leq C\int^\tau_0h(t)\mathcal{E}(\rho_{\epsilon},\mathbf{u}_{\epsilon}|r,\mathbf U)dt+\delta\int^\tau_0\|\nabla_x\mathbf v_{\epsilon}-\nabla_x\mathbf V\|^2_{L^{2}(\Omega)}
+\|\partial_z\mathbf v_{\epsilon}-\partial_z\mathbf V\|^2_{L^{2}(\Omega)}dt\nonumber\\
&\hspace{8pt}+\int^\tau_0\int_\Omega\frac{\rho_{\epsilon}}{r}(\mathbf V-\mathbf v_{\epsilon})(\Delta_x\mathbf V+\partial_{zz}\mathbf V)dxdzdt
-\int^\tau_0\int_\Omega\frac{\rho_{\epsilon}}{r}(\mathbf V-\mathbf v_{\epsilon})\nabla_xp(r)dxdzdt\nonumber\\
&\hspace{8pt}-\int^\tau_0\int_{\Omega}P''(r)((\rho_{\epsilon}-r)\partial_tr+\rho_{\epsilon}\mathbf v_{\epsilon}\nabla_xr)dxdzdt
-\int^\tau_0\int_\Omega p(\rho_{\epsilon})\text{div}_x\mathbf Vdxdzdt+o(\epsilon^2).
\end{align*}

Then we deduce that
\begin{align}
\mathcal{E}&(\rho_{\epsilon},\mathbf{u}_{\epsilon}|r,\mathbf U)|^{t=\tau}_{t=0}+\int^\tau_0\int_\Omega\big{(}(\nabla_x\mathbf v_{\epsilon}-\nabla_x\mathbf V):(\nabla_x\mathbf v_{\epsilon}-\nabla_x\mathbf V)
+|\partial_z\mathbf v_\epsilon-\partial_z\mathbf V_\epsilon|^2+\epsilon^2|\nabla w_{\epsilon}|^2\big{)}dxdzdt\nonumber\\
&\leq C\int^\tau_0h(t)\mathcal{E}(\rho_{\epsilon},\mathbf{u}_{\epsilon}|r,\mathbf U)dt+\delta\int^\tau_0\|\nabla_x\mathbf v_{\epsilon}-\nabla_x\mathbf V\|^2_{L^{2}(\Omega)}
+\|\partial_z\mathbf v_{\epsilon}-\partial_z\mathbf V\|^2_{L^{2}(\Omega)}dt\nonumber\\
&\hspace{8pt}+\int^\tau_0\int_\Omega(\frac{\rho_{\epsilon}}{r}-1)(\mathbf V-\mathbf v_{\epsilon})(\Delta_x\mathbf V+\partial_{zz}\mathbf V)dxdzdt
-\int^\tau_0\int_\Omega\frac{\rho_{\epsilon}}{r}(\mathbf V-\mathbf v_{\epsilon})\nabla_xp(r)dxdzdt\nonumber\\
&\hspace{8pt}-\int^\tau_0\int_{\Omega}P''(r)((\rho_{\epsilon}-r)\partial_tr+\rho_{\epsilon}\mathbf v_{\epsilon}\nabla_xr)dxdzdt
-\int^\tau_0\int_\Omega p(\rho_{\epsilon})\text{div}_x\mathbf Vdxdzdt+o(\epsilon^2).\label{4.12}
\end{align}

\subsection{Step 3}
We are now in a position to estimate the remaining terms in the relative energy inequality (\ref{4.12}). It is clear to check that
\begin{align}
-\int^\tau_0&\int_\Omega\frac{\rho_{\epsilon}}{r}(\mathbf V-\mathbf v_{\epsilon})\nabla_xp(r)+p(\rho_{\epsilon})\text{div}_x\mathbf V+P''(r)((\rho_{\epsilon}-r)\partial_tr+\rho_{\epsilon}\mathbf v_{\epsilon}\nabla_xr)dxdzdt\nonumber\\
&=-\int^\tau_0\int_\Omega(\rho_{\epsilon}-r)P''(r)\partial_tr+P''(r)\rho_{\epsilon}\mathbf v_{\epsilon}\cdot\nabla_xr+\rho_{\epsilon} P''(r)(\mathbf V-\mathbf v_{\epsilon})\cdot\nabla_xr+p(\rho_{\epsilon})\text{div}_x\mathbf Vdxdzdt\nonumber\\
&=-\int^\tau_0\int_\Omega(\rho_{\epsilon}-r)P''(r)\partial_tr+P''(r)\rho_{\epsilon}\mathbf V\cdot\nabla_xr+p(\rho_{\epsilon})\text{div}_x\mathbf Vdxdzdt\nonumber\\
&=-\int^\tau_0\int_\Omega\rho_{\epsilon} P''(r)(\partial_tr+\mathbf V\cdot\nabla_xr)
-rP''(r)\partial_tr+p(\rho_{\epsilon})\text{div}_x\mathbf Vdxdzdt\nonumber\\
&=-\int^\tau_0\int_\Omega\rho_{\epsilon} P''(r)(-r\text{div}_x\mathbf V-r\partial_zW)
-rP''(r)\partial_tr+p(\rho_{\epsilon})\text{div}_x\mathbf Vdxdzdt\nonumber\\
&=-\int^\tau_0\int_\Omega\text{div}_x\mathbf V\big{(}p(\rho_{\epsilon})-p'(r)(\rho_{\epsilon}-r)-p(r)\big{)}dxdzdt
+\int^\tau_0\int_\Omega p'(r)(\rho_{\epsilon}-r)\partial_zWdxdzdt,
\end{align}
where we have used the fact that $\partial_tr+\text{div}_x\mathbf Vr+\mathbf V\cdot\nabla_xr+r\partial_zW=0$.

Using the analogous argument as in \cite{ma} Section 2.2.5, we can easily carry out the following estimate:
\begin{align}
|\int^\tau_0\int_\Omega\text{div}_x\mathbf V\big{(}p(\rho_{\epsilon})-p'(r)(\rho_{\epsilon}-r)-p(r)\big{)}dxdzdt|
\leq C\int^\tau_0h(t)\mathcal{E}(\rho_{\epsilon},\mathbf u_{\epsilon}|r,\mathbf U)dt.
\end{align}

According to the periodic boundary condition, it follows that
\begin{align}
\int^\tau_0\int_\Omega p'(r)(\rho_{\epsilon}-r)\partial_zWdxdzdt
=\int^\tau_0dt\int_{\mathbb{T}^2}(\int^1_0\partial_zWdz)p'(r)(\rho_{\epsilon}-r)dx=0.
\end{align}

Furthermore, an argument similar to the one used in \cite{kr} Section 6.3 shows that
\begin{align}
\int_\Omega&(\frac{\rho_{\epsilon}}{r}-1)(\mathbf V-\mathbf v_{\epsilon})(\Delta_x\mathbf V+\partial_{zz}\mathbf V)dxdz\nonumber\\
&\leq C\mathcal{E}(\rho_{\epsilon},\mathbf u_{\epsilon}|r,\mathbf U)+\delta\|\nabla_x\mathbf v_{\epsilon}-\nabla_x\mathbf V\|^2_{L^2}
+\delta\|\partial_z\mathbf v_{\epsilon}-\partial_z\mathbf V\|^2_{L^2}.
\end{align}

Therefore, putting $(4.12)-(4.16)$ together, we have
\begin{align}
\mathcal{E}(\rho_{\epsilon},\mathbf u_{\epsilon}|r,\mathbf U)(\tau)\leq C\int^\tau_0h(t)\mathcal{E}(\rho_{\epsilon},\mathbf u_{\epsilon}|r,\mathbf U)(t)dt+o(\epsilon^2).
\end{align}

Then applying the Gronwall's inequality, we finish the proof of Theorem 3.1.

\vskip 0.5cm


\vskip 0.5cm
\noindent {\bf Acknowledgements}

\vskip 0.1cm
The research of H. G is partially supported by the NSFC Grant  No. 11531006. The research of \v S.N. is supported by the Czech Sciences Foundation (GA\v CR),   GA19-04243S and RVO 67985840.  The research of T.T. is supported by the NSFC Grant No. 11801138. The authors thank Prof. Edriss Titi's wonderful and great talk on the course"compact course Mathematical Analysis of Geophysical Models and Data Assimilation" held from 26 June to 10 July 2020, which offers insightful and constructive suggestions. We would like to thank to Prof. A. Novotn\' y for his remarks and suggestions.


\end{document}